\documentclass[11pt]{article}
\usepackage{epsfig}
\usepackage{t1enc}
    \usepackage[latin1]{inputenc}
    \usepackage[english]{babel}
        \usepackage{latexsym}
\usepackage{amssymb}
\usepackage{euscript}

\def\begg{\begin{equation}}
\def\endd{\end{equation}}

\newcommand{\ep}{{\epsilon}}
\begin{document}

\setlength{\arraycolsep}{.136889em}
\renewcommand{\theequation}{\thesection.\arabic{equation}}
\newtheorem{thm}{Theorem}[section]
\newtheorem{propo}{Proposition}[section]
\newtheorem{lemma}{Lemma}[section]
\newtheorem{corollary}{Corollary}[section]
\newtheorem{remark}{Remark}[section]
\newtheorem{consequence}{Consequence}[section]

\smallskip

\centerline{\Large\bf ABOUT THE  DISTANCE }
\smallskip
\centerline{\Large\bf BETWEEN RANDOM WALKERS ON SOME GRAPHS}

\bigskip\bigskip
\bigskip\bigskip

\bigskip\bigskip
\renewcommand{\thefootnote}{1}

\noindent
\textbf{Endre Cs\'{a}ki}\footnote{ Research supported by the Hungarian National
Research, Development and Innovation Office - NKFIH K 108615. }
\newline
Alfr\'ed R\'enyi Institute of Mathematics, Hungarian
Academy of Sciences,
\newline
Budapest, P.O.B. 127, H-1364, Hungary. E-mail address:
csaki.endre@renyi.mta.hu

\bigskip
\renewcommand{\thefootnote}{2}
\noindent
{\textbf{Ant\'{o}nia
F\"{o}ldes}\footnote{Research supported by a PSC CUNY Grant, No.
69040-0047.}}
\newline
Department of Mathematics, College of Staten
Island, CUNY, 2800 Victory Blvd., Staten Island, New York 10314,
U.S.A.  E-mail address: antonia.foldes@csi.cuny.edu

\bigskip
\renewcommand{\thefootnote}{3}
\noindent
\noindent
\textbf{P\'al R\'ev\'esz}\footnote{ Research supported by the Hungarian National
Research, Development and Innovation Office - NKFIH K 108615. }
\newline
\noindent
Institut f\"ur
Statistik und Wahrscheinlichkeitstheorie, Technische Universit\"at
Wien, Wiedner Hauptstrasse 8-10/107 A-1040 Vienna, Austria.
E-mail address: revesz@ci.tuwien.ac.at
\bigskip

\medskip
\noindent{\bf Abstract}
\newline
We consider two or more simple symmetric walks on  $\mathbb{Z}^d$ 
and the 2-dimensional comb lattice, and investigate the
properties of the distance among the walkers.
\medskip


\medskip

\vspace{.1cm}

\medskip

\section{Introduction }
\renewcommand{\thesection}{\arabic{section}} \setcounter{equation}{0}
\setcounter{thm}{0} \setcounter{lemma}{0}

Almost a hundred years ago, P\'olya \cite{P} in 1921 proved that
on  $\mathbb{Z}^1$ and  $\mathbb{Z}^2$ simple random walks are recurrent and  
two independent walkers meet infinitely often with probability one, but on 
$\mathbb{Z}^d$, for $d \geq 3$, simple random walks are transient 
and two independent walkers meet only finitely often with probability one. 
Nowadays meeting at  the same place at the same time is called a collision, 
so we will use this term, to avoid any confusion. On $\mathbb{Z}$ two 
walkers not only collide infinitely many times, but they collide even in the 
origin infinitely many times with probability one. In their landmark paper 
Dvoretzky and Erd\H{o}s \cite{DE} in 1950  recall the celebrated P\'olya
result and, among their final remarks, they mention without proof that on
$\mathbb{Z}$ three independent random walkers collide (all three together)
infinitely often with probability one. A short elegant proof was given for
this statement in Barlow, Peres and Sousi \cite{BPS} in 2012. However four
walkers in $\mathbb{Z},$  three walkers in $\mathbb{Z}^2$ and  two walkers
 in $\mathbb{Z}^d,$ with $d\geq 3$ will only collide
finitely many times with probability 1. Khrishnapur and Peres \cite{KP} in
2004 studied this problem on the comb lattice. They proved that even though the
comb is recurrent, two independent random walkers on the comb lattice collide
only finitely often with probability 1.

Our question is that in case the walkers collide only finitely often, then how
does their distance grow as a function of time. We want to establish upper
class results for the distance of two or more walkers and lower class results
for the distance of four or more walkers on $\mathbb{Z}$. Similarly, for the
distance in $\mathbb{Z}^2,$ we give upper class results for the distance of  
two or more walkers  and lower class results for the distance of three or more
walkers. For $\mathbb{Z}^d,$ if $d\geq 3,$ we get upper and lower class 
results for  the distance of two or more  walkers. Finally, we will 
investigate the distance of two or more walkers on the comb lattice.

We start with some definitions. Let $\bf{G}$  be a connected graph with
vertex set $\bf{V(G)}.$ Two neighboring connected vertices $v$ and $w$ form
an edge of $\bf{G}.$ A random walk $S(n)$  on $\bf{G}$ is defined with the
following one step transition probabilities
\begg  p(u,v):=P(S(n+1)=v|S(n)=u)=\frac{1}{deg(u)}  \label{nocsak} \endd
for neighboring vertices  $u$ and $v$ in $\bf{V(G)},$ where $deg(u)$ 
is the number of neighbours of $u,$ otherwise $p(u,v)=0.$
We define the {\it graph distance}, which we will simply call {\it distance}, 
of $u$ and $v$ in $\bf{V(G)}$ as the minimal number of steps the walker needs 
to arrive from $u$ to $v.$ Formally,
\begg
d(u,v):=\min\{k>0:\, P(S(n+k)=v|S(n)=u)>0\}.  \label{dist}
\endd
Or, equivalently, the distance $d(u,v)$ is the length of the shortest path from
vertex $u$ to vertex $v$  in $\bf{V(G)}.$ In $\mathbb{Z}^d$ we will use
Euclidean distance. In $\mathbb{Z}$ however these two distances are the same.

\section{Preliminary results}
\renewcommand{\thesection}{\arabic{section}} \setcounter{equation}{0}
\setcounter{thm}{0} \setcounter{lemma}{0}
In this section we list some known important results which we will need later
on. Put ${\bf W}(t)=(W_1(t), W_2(t),...,W_d(t)),$  where
$W_1(t), W_2(t),...,W_d(t)$ are independent standard Wiener processes. Then the
${\mathbb R}^d$  valued process ${\bf W}(t)$ is called the standard
$d$-dimensional Wiener process. Let ${\bf S}(n)$ be the location of a walker in
${\mathbb Z}^d$ after $n$-steps, where  the simple symmetric walk
${\bf S}(n)={\bf X}(1)+{\bf X}(2)+...+{\bf X}(n),$ $n=1,2,...$ and
${\bf X}(1), {\bf X}(2), ...,{\bf X}(n)$ are i.i.d. random vectors with
${\bf S}(0)=0$,
\begg
P({\bf X}(1)=e_i) =P({\bf X}(1)= -e_i)=\frac{1}{2d}, \,\, i=1,2,...,d,
\label{vector}
\endd
where $e_1, e_2,...,e_d$ are the orthogonal unit-vectors  in ${\mathbb Z}^d.$

In Section 3 we investigate the distance on ${\mathbb Z}^d$, and show that
the upper bound is a consequence of the law of the iterated logarithm (LIL),
and that the lower bound can be established using a result of
Dvoretzky and Erd\H os \cite{DE}.
 For the multidimensional LIL we refer, e.g., to R\'ev\'esz \cite{R13},
Theorem 19.1.

\noindent
{\bf Theorem A} {\it For the  $d$-dimensional  standard Wiener process and the 
simple symmetric random walk we have for any $d \geq 1$
\begg
\limsup_{t\to\infty}\frac{\|{\bf W}(t)\|}{\sqrt{t\log\log t}}=
\sqrt{2}\quad a.s.
\endd
and}
\begg
\limsup_{n\to\infty}\frac{\|{\bf S}(n)\|}{\sqrt{n\log\log n}}=
\sqrt{\frac{2}{d}}\quad a.s.
\endd

We use the following definition (cf. R\'ev\'esz \cite{R13}, page 36): The 
function $g(t)$ belongs to the lower-lower class (LLC) of the random process
$\{Y(t),\, t\geq 0\}$ if for almost all $\omega\in \Omega$ there exists  a
$t_0=t_0(\omega)$ such that $Y(t)>g(t)$ if $t>t_0$.
\smallskip

\noindent
{\bf Theorem B} Dvoretzky-Erd\H{o}s \cite{DE} {\it Let $a(t)$ be a
nonincreasing nonnegative function. Then, for the $d$-dimensional random walk 
${\bf S}(n)$ and the standard Wiener process ${\bf W}(t)$
\begin{eqnarray}
t^{1/2}a(t)&\in& {\rm LLC}(\|{\bf W}(t)\|)\qquad  (d \geq 3)\\
n^{1/2}a(n)&\in& {\rm LLC}(\|{\bf S}(n)\|)\qquad  (d \geq 3)
\end{eqnarray}
if and only if }
$$\sum_{n=1}^{\infty} (a(2^n))^{d-2}<\infty.$$
Here and throughout $\|\cdot\|$ denotes the Euclidean distance in $d$ 
dimensions.
\begin{remark}
Note that the same results hold when $a(n)$ is replaced by $c\,a(n)$
with any positive constant $c$.
\end{remark}

We will need a special case of the multidimensional invariance principle,
explicitly stated in R\'ev\'esz \cite{R13}, Theorem 18.2. (There are more 
precise results in the literature, but we don't need them here.)

\smallskip
\noindent
{\bf Theorem C} {\it On a rich enough  probability space
one can define a standard $d$-dimensional Wiener process
$\{{\bf W}(t), \, t>0\}$ and a simple random walk
$\{{\bf S}(n),\,n=0,1,2...\}$ on ${\mathbb Z}^d$ such that} for $d\geq 1$
$$
\|{\bf S}(n)-{\bf W}(n/d)\|=O(n^{1/4}(\log\log n)^{3/4})\quad {\rm a.s.}
\quad {\it as}\,\,\, n\, \to \infty.
$$

\smallskip
We will use the following result from linear algebra.

\smallskip \noindent
{\bf Lemma D} \emph{ Let $\mathbb{S}$ be a vector space with basis
$\{w_i,\,\, i=1,2,...,n\}.$ Define two subspaces} $\mathbb{U}$  \emph{and}
$\mathbb{V}$ \emph{in}
$\mathbb{S}$ \emph{by}
\begin{eqnarray}
\mathbb{U}&=&Span\{w_1-w_i,\,\, i=2,3,...,n\} \nonumber \\
\mathbb{V}&=&Span\{ \sum_{j=1}^i (w_j-w_{i+1}),\,\,i=1,2,...,n-1\}.
\end{eqnarray}
\emph{Then the two subspaces defined above are the same.}

The results presented so far are needed for proving our theorems about the
$d$-dimensional random walk in Section 3. In Section 4 we consider the case
of simple random walk on  the 2-dimensional comb.

The 2-dimensional comb lattice
${\mathbb C}^2$  is obtained from ${\mathbb Z}^2$ by removing all
horizontal edges off the $x$-axis. In this context the $x$-axis is usually
called the backbone of the comb and the vertical lines are called teeth.
A formal way of describing a simple random walk ${\bf C}(n)$ on the above
2-dimensional comb lattice ${\mathbb C}^2$ can be formulated via its
transition probabilities as follows: for any integers $x$ and $y$ define
\begin{equation}
 P({\bf C}(n+1)=(x,y\pm 1)\mid {\bf C}(n)=(x,y))=\frac12,
\quad {\rm if}\,\,  y\neq 0,
\end{equation}
\begin{equation}
 P({\bf C}(n+1)=(x\pm 1,0)\mid {\bf C}(n)=(x,0))=
 P({\bf C}(n+1)=(x,\pm 1)\mid {\bf C}(n)=(x,0))=\frac14.
\end{equation}

A  compact  way of describing the just introduced transition probabilities for
this simple random walk ${\bf C}(n)$ on ${\mathbb C}^2$ is via (\ref{nocsak}).

As far as we  know, the first discussion of random walk  on  the comb was
given by Weiss and Havlin \cite{WH}. Bertacchi and Zucca \cite{BZ0}
obtained the following space-time asymptotic estimates for the $n$-step
transition probabilities $p^{{\mathbb C}^2}(u,v,n)$, where $u=(x_1,y_1)$
and $v=(x_2,y_2)$ are two vertices on the comb.
For any $u$ and $v$ fixed vertices on ${\mathbb C}^2,$ with $deg(\cdot)$ as  
in (\ref{nocsak}),
\begg
p^{{\mathbb C}^2}(u,v,n)\sim
\frac{2^{\frac{1}{4}-1} deg(v)}{\Gamma(\frac{1}{4})n^{\frac{3}{4}}},\qquad
{\rm as} \,\, n \to \infty \label{berzu}
\endd
whenever $n+d(u,v)$ is even, and $p^{{\mathbb C}^2}(u,v,n)=0$ if
$n+d(u,v)$ is odd,  and $\sim$ stands for asymptotic equality.

Here we recall our construction  of  two dimensional comb walk  from
\cite{CCFR08}, which we used there to  prove Theorem K below, as we will need 
some parts of this construction later on. Consider a sample space large 
enough to contain  two independent simple symmetric random walks  $S_1(n)$ 
and  $S_2(n),\, n=1,2,\ldots,$ on the integer lattice on the line, and an 
i.i.d. sequence of geometric random variables $G_i, i=1,2,\ldots,$ with

$$ P(G_1=k)=\frac1{2^{k+1}},\quad k=0,1,2,\ldots
$$ which is independent from the two random walks as well.
We may then construct a simple random walk on
the 2-dimensional comb lattice ${\mathbb C}^2$ as follows. Let $\rho_2(N)$
be the time of the $N$-th return  to zero of the second random walk
$S_2(\cdot)$, i.e., $\rho(0)=0,$ and

$$\rho_2(N):= \min\{j>\rho_2(N-1) : S_2(j)=0\}.$$
\noindent Put
$T_N=G_1+G_2+\ldots G_N$,
$N=1,2,\ldots\,\,.$ For
$n=0,\ldots, T_1$, let $C_1(n)=S_1(n)$ and
$C_2(n)=0$. For $n=T_1+1,\ldots, T_1+\rho_2(1)$, let $C_1(n)=C_1(T_1)$,
$C_2(n)=S_2(n-T_1)$. In general, for $T_N+\rho_2(N)<n\leq
T_{N+1}+\rho_2(N)$, let
$$
C_1(n)=S_1(n-\rho_2(N)),
$$
$$
C_2(n)=0,
$$
and, for $T_{N+1}+\rho_2(N)<n\leq T_{N+1}+\rho_2(N+1)$, let
$$
C_1(n)=C_1(T_{N+1}+\rho_2(N))=S_1(T_{N+1}),
$$
$$
C_2(n)=S_2(n-T_{N+1}).
$$
Then it can be seen in terms of these definitions for $C_1(n)$ and
$C_2(n)$ that ${\bf C}(n)=(C_1(n),C_2(n))$ is a simple random
walk on the 2-dimensional comb lattice ${\mathbb C}^2$. Define the local time
  of $S_i(n), \,\,i=1,2, $ at zero by
$$
\xi_i(0,n):=\sum_{j=1}^n I\{S_i(j)=0\}
$$
and denote the
number of horizontal and vertical steps of ${\bf C}(n)$ by $H_n$ and $V_n$
respectively, with $H_n+V_n=n.$ Clearly $H_n$ is the sum of $\xi_2(0,V_n) $
i.i.d. geometric random variables $G_i$ as described above, where the last
geometric random variable may be truncated.
Moreover  ${\bf C}(n)=(C_1(n),C_2(n))=(S_1(H_n),S_2(V_n)).$

We will also need the following increment result  of
Cs\"org\H{o} and R\'ev\'esz (\cite{CR81}, page 115) for a random walk.

\smallskip
\noindent
{\bf Theorem E}  {\it Let  $X_1, X_2, \cdots$ be a sequence of i.i.d. random
variables with mean zero and variance one, satisfying the following
condition:
there exists a $t_0$ such that $E(e^{t X_1})$ is finite for $|t| < t_0. $

\noindent
Let $0< a_N$ be a non-decreasing sequence of integers such that
$N/a_N$  is also non-decreasing and $N/a_N \to \infty$. Then, for
$S(n)=X_1+X_2 +\cdots +X_n$, as $N\to \infty,$ we have almost surely that}

$$\max_{0\leq N-a_N} \max_{k\leq a_N}|S(n+k)-S(n)|=O(a_N^{1/2}(\log(N/a_N)
+\log\log N)^{1/2}).$$

The following theorem is a version of
Hoeffding's inequality, which is stated explicitly  in \cite{TO}.

\smallskip \noindent
{\bf Theorem F} {\it Let $G_i$ be i.i.d.  random variables with the
common geometric distribution $ P(G_i=k)=2^{-k-1},\quad k=0,1,2,...\,\,.$
Then

$$ P\left(\max_{1\leq j\leq n}
\left|\sum_{i=1}^j (G_i-1)\right|>\lambda\right)\leq
2\exp(-\lambda^2/{8n})
$$
for $0<\lambda<na,$ with some} $a>0.$

\smallskip \noindent
\begin{consequence}  For any small $0<\delta< a$  and $n$ big enough
$$
 P\left(\sum_{i=1}^n G_i\geq (1+\delta)n \right)\leq \frac{1}{n^2}.
$$
\end{consequence}

Let $S(n)$ be a simple symmetric random walk on the line with local
time $\xi(0,n)$ at zero. 
We recall from Cs\'aki and F\"oldes \cite{CsF83} the following result.

\smallskip \noindent
{\bf Theorem G} {\it Suppose that $x_n \to \infty,$
\,$\displaystyle{\frac{x_n}{n^{1/2}} \to 0}$ as $ n \to \infty.$
Then, for any $0<\ep\leq 1,$
$$
P(\xi(0,n)\geq x_n n^{1/2}) \leq c \exp\left(- \frac{(1-\ep)x^2_n}{2}\right)$$
for some constant $c$  if $n$ is large enough.}
\noindent
\begin{remark} In what follows, we will disregard parity issues. Namely
we will use results for the $n-$step transition probability for arbitrary
$n,u$ and $v,$ when these results are only proved for even $n,u$ and $v.$
When these values are not even, and the corresponding probability is not zero,
then their asymptotic behavior is the same as for even values (see Bertacchi
and Zucca {\rm{\cite {BZ}, Sections 3, 4 and 10}}).
\end{remark}

\smallskip
For completeness, we present the following trivial lemma, which is likely
well-known.

\smallskip \noindent
{\bf Lemma H} {\it Let $u$ and $v$ be two distinct vertices on ${\mathbb C}^2.$
For any $n$ such that $n+d(u,v)$ is even}
\begg
deg(u) p^{{\mathbb C}^2}(u,v,n)= deg(v)p^{{\mathbb C}^2}(v,u,n),
\endd
{\it where $deg(\cdot)$ is as in} (\ref{nocsak}).

\noindent
{\bf Proof.} Denote the set of all $n$-steps paths connecting vertex $u$
to vertex $v$ through the vertices $z_1,z_2,...,z_{n-1}$ by
$\lambda_i:=(u,z^i_1,z^i_2,...z^i_{n-1},v),$ and the set of all paths of length
$n$ from $u$ to $v$ by $\Lambda(u,v,n)$ and, for any vertex $z,$
let $q(z)=\frac{1}{deg(z)}.$ Then
$$p^{{\mathbb C}^2}(u,v,n)=\sum_{\lambda_i\in \Lambda(u,v,n)}
q(u)q(z^i_1)q(z^i_2)...q(z^i_{n-1}),$$
$$p^{{\mathbb C}^2}(v,u,n)=\sum_{\lambda^*_i\in \Lambda(v,u,n)}
q(v) q(z^i_{n-1}) q(z^i_{n-2})...q(z^i_1), $$
where $\lambda_i^*$ is the reversed path of $\lambda_i.$
Hence, the above two probabilities are equal if $q(u)=q(v),$ and  differ in a
factor 2  if  $u$ and  $v$ are having different degrees. $\Box$

In (\ref{berzu}) the two vertices are fixed. We will need a more general
estimate for the $n$-step transition probabilities
$p^{{\mathbb C}^2}((0,0), (0,k),n)$ .  Define
$$\kappa=k/n \,\,\,{\rm and}\,\,\,\, \phi(\kappa)
=\log((1-\kappa)^{\kappa-1}(1+\kappa)^{-\kappa-1}).$$

\noindent
We recall the second statement from Theorem 5.5 of Bertacchi and Zucca
\cite{BZ}.

\smallskip \noindent
{\bf Theorem I} { \it If $\kappa \in [0,n^{-1/2-\ep}]$ for some $\ep>0$, then
as} $n\to \infty$
$$p^{{\mathbb C}^2}((0,2k),(0,0),2n)\sim
\frac{\sqrt{2}e^{n\phi(\kappa)}}{\Gamma(1/4)n^{3/4}},$$
{\it uniformly with respect to} $\kappa \in [0,n^{-1/2-\ep}]$.

A simple calculation shows that $\phi(\kappa)$ is decreasing  for
$\kappa\geq 0$, hence $\phi(\kappa)\leq \phi(0)=0.$
Consequently $e^{n\phi(\kappa)}\leq 1$.  Combining this with Lemma H,
the following obtains.

\smallskip \noindent
\begin{consequence} For $k\leq n^{1/2-\ep}$ and some c>0,
$$p^{{\mathbb C}^2}((0,0),(0,k),n)\leq \frac{c}{n^{3/4}}.$$
\end{consequence}

For any two vertices $u$ and $v$ on ${\mathbb C}^2,$ we define the Green
function associated with the random walk on the comb as
$$G(u,v|z):=\sum_{n= 0}^{\infty} p^{{\mathbb C}^2}(u,v,n)z^n.$$
Bertacchi and Zucca  in \cite{BZ} show that for $u=(0,0),\,$
$G((0,0),(k,\ell)|z)$ can be given explicitly as follows:

\[G((0,0),(k,\ell)|z)=\left \{ \begin{array}{ll}
\frac{1}{2}G(z)(F_1(z))^{|k|}(F_2(z))^{|\ell|}\,\, &\mbox {if $\ell\neq 0$ }
\\
G(z)(F_1(z))^{|k|}\,\, &\mbox {if $\ell=0$},
\end{array}
\right. \]
where
\begin{eqnarray*}
G(z)&=&\frac{\sqrt 2}{\sqrt{1-z^2+\sqrt{1-z^2}}},\nonumber\\
F_1(z)&=& \frac{1+\sqrt{1-z^2}-\sqrt{2}\sqrt{1-z^2+\sqrt{1-z^2}}}{z},\nonumber\\
F_2(z)&=& \frac{1-\sqrt{1-z^2}}{z}.
\end{eqnarray*}
We will use this elegant result to get the asymptotic behavior of the
probability $P(C_2(n)=0).$ Selecting $\ell=0$ and  summing for all
$k=0,\pm 1,\pm2,...,$ we easily obtain the generating function of
$P(C_2(n)=0).\,$ Namely, with the notation
 $$H(z):=\sum_{n= 0}^{\infty} P(C_2(n)=0)z^n,$$
we have
$$H(z)=\sum_{k=-\infty}^{\infty} G((0,0)(k,0)|z)
=G(z)+2\sum_{k=1}^{\infty}G(z)(F_1(z))^k=G(z)\frac{1+F_1(z)}{1-F_1(z)}.$$
Now, just like as it is used in \cite{BE}, we may  apply the
Hardy-Littlewood-Karamata theorem in the following form: If
$$ F(z)=\sum a_n z^n \sim \frac{C}{(1-z)^{\alpha}} \quad {\rm as}
\quad z\to 1^- \quad {\rm with}\quad \alpha \notin \{0,-1,-2,...\} ,$$
and if $F(z)$ is analytic in some domain with the exception of $z=1,$ then
$$a_n\sim \frac{C}{\Gamma(\alpha)}n^{\alpha-1}\quad {\rm as}
\quad n\to \infty.$$
An easy calculation yields the asymptotic behavior of $H(z)$. We only
give a short indication of this calculation:
\begin{eqnarray*}
G(z)\frac{(1+F_1(z))}{1-F_1(z)} &=&
\frac{\sqrt 2}{\sqrt{1-z^2+\sqrt{1-z^2}}}\,
\frac{z+1+\sqrt{1-z^2}-\sqrt{2}\sqrt{1-z^2+\sqrt{1-z^2}} }
{z-1-\sqrt{1-z^2}+\sqrt{2}\sqrt{1-z^2+\sqrt{1-z^2}} }\\
&\sim &\frac{\sqrt 2}{(1-z^2)^{1/4}}\,\frac{2}{\sqrt 2(1-z^2)^{1/4}} \\
&=&\frac{2}{\sqrt{(1-z)(1+z)}}\\
&\sim& \frac{\sqrt{2}}{\sqrt{1-z}}, \quad {\rm as}\quad z\to 1^-. \\
\end{eqnarray*}
\noindent
Thus, we have $\alpha=1/2,$ implying the following result.

\noindent
\begin{consequence} For the second coordinate  $C_2(n)$  of the comb walk we 
have
$$P(C_2(n)=0)\sim \frac{\sqrt{2}}{\sqrt{\pi n}}, \quad{\rm as}\,\,\,
n\to \infty.$$
\end{consequence}

A further insight to the nature of the random walk on a comb was
provided by Bertacchi \cite{BE}, who established the following remarkable
weak convergence result for the walk
${\bf C}(n)=(C_1(n),C_2(n))$ on the comb ${\mathbb C}^2$.

\medskip
\noindent
{\bf Theorem J} {\it For the random walk $\{{\bf C}(n)=(C_1(n),C_2(n));
n=0,1,2,\ldots\}$ on ${\mathbb C}^2$,we have}
\begin{equation}
\left(\frac{C_1(nt)}{n^{1/4}},\frac{C_2(nt)}{n^{1/2}};\, t\geq 0\right)
{\buildrel{\rm Law}\over\longrightarrow}\, (W_1(\eta_2(0,t)), W_2(t);\,
t\geq 0), \quad n\to\infty,
\label{ber}
\end{equation}
{\it where $W_1$, $W_2$ are two independent standard Wiener processes (Brownian
motions) and $\eta_2(0,t)$ is the local time process of $W_2$ at zero,
and ${\buildrel{\rm Law}\over\longrightarrow}$ denotes weak convergence
on $C([0,\infty),{\mathbb R}^2)$ endowed with the topology
of uniform convergence on compact intervals.}

For the definition of $\eta_2(0,t)$ see e.g. R\'ev\'esz \cite{R13},  page 107.

In our paper \cite{CCFR08} we gave a joint strong
approximation result for the two coordinates of this walk.

\medskip
\noindent
{\bf Theorem K}
{\it On an appropriate probability space for the random walk
$\{{\bf C}(n)=(C_1(n),C_2(n));\newline n=0,1,2,\ldots\}$ on ${\mathbb
C}^2$, one can construct two independent standard Wiener processes
$\{W_1(t);\, t\geq 0\}$, $\{W_2(t);\, t\geq 0\}$ so that, as $n\to\infty$,
we have with any $\varepsilon>0$
$$
n^{-1/4}|C_1(n)-W_1(\eta_2(0,n))|+n^{-1/2}|C_2(n)-W_2(n)|
=O(n^{-1/8+\varepsilon})\quad {\rm a.s.},
$$
where $\eta_2(0,\cdot)$ is the local time process at zero of
$W_2(\cdot)$.}

From the many consequences of this result, we will need the following two.

\begin{corollary} For the horizontal and vertical coordinates of
${\bf C}(n)=(C_1(n),C_2(n))$ we have
\begin{equation}
\limsup_{n\to\infty}\frac{|C_1(n)|}{n^{1/4}(\log\log n)^{3/4}}
=\frac{2^{5/4}}{3^{3/4}}\quad {\rm a.s.},
\label{itlog1}
\end{equation}
\begin{equation}
\limsup_{n\to\infty}\frac{|C_2(n)|}{(2n\log\log n)^{1/2}}=1 \quad
{\rm a.s.}
\label{itlog2}
\end{equation}
\end{corollary}

\smallskip
For the distribution of the hitting time of a simple random walk, we need the
following result (cf., e.g., Feller \cite{FE}, Ch. 3.7, Theorem 2 and
Theorem 3).

\medskip
\noindent
{\bf Theorem L}
{\it Let $\{S(i),\, i=0,1,2,\ldots\}$ be a simple symmetric random walk
on the line with $S(0)=0$, and define the hitting time}
\begg
\beta(r)=\min\{i>0:\, S(i)=r\},
\label{beta}
\endd
{\it where $r$ is a positive integer. Then}
\begg
P(\beta(r)=N)=\frac{r}{N}{N\choose\frac{N+r}{2}}2^{-N},\quad N=r,r+1,\ldots
\label{feller1}
\endd
{\it whith $N+r$ even, and
 \begg
\lim_{r\to\infty}P(\beta(r)<ur^2)=\sqrt{\frac{2}{\pi}}
\int_{1/\sqrt{u}}^\infty e^{-s^2/2}\, ds, \quad u>0.
\label{feller2}
\endd}

\section{Distance on $\mathbb{Z}^d$}
\renewcommand{\thesection}{\arabic{section}} \setcounter{equation}{0}
\setcounter{thm}{0} \setcounter{lemma}{0}

Let $\{{\bf S}_i(\cdot),\quad i=1,2,...,K\},$ be $K$ independent random walks 
on $\mathbb{Z}^d,$ the paths of the $K$ walkers.
We consider the maximal distance between $K$ walkers as follows.
$$D_K^{\mathbb{Z}^d}(n):=\max_{i\neq j,\,\, i,j\leq K}
\|{\bf S}_i(n)-{\bf S}_j(n)\|,$$
where $\|\cdot\|$ denotes Euclidean distance.

Similarly, for $\{{\bf W}_i(\cdot),\quad i=1,2,...,K\},$ $K$ independent 
standard $d$-dimensional Wiener processes, all starting from $0$, let
$$D_K^{\mathbb{R}^d}(t)=
\max_{i\neq j,\,\, i,j\leq K} \|{\bf W}_i(t)-{\bf W}_j(t)\|.$$

Concerning  upper class results,  we prove our next result from the law of 
the iterated logarithm (LIL).

\begin{thm}
For $K\geq 2$ and $d=1,2,\ldots,$ we have
\begg
\limsup_{t\to\infty}\frac{D_K^{\mathbb{R}^d}(t)}{\sqrt{t\log\log t}}
=2\quad a.s.
\label{iterated1}
\endd
and
\begg
\limsup_{n\to\infty}\frac{D_K^{\mathbb{Z}^d}(n)}{\sqrt{n\log\log n}}
=\frac{2}{\sqrt{d}}\quad a.s.
\label{iterated2}
\endd
\end{thm}

For the lower classes, we prove the following results.
\begin{thm}
For $d=1,2,\ldots$ and $K\geq 1+\frac{3}{d}$ we have
\begg
\sqrt{t} a(t) \in {\rm LLC}(D_K^{\mathbb{R}^d}(t))
\label{lower1}
\endd
and
\begg
\sqrt{n} a(n) \in {\rm LLC}(D_K^{\mathbb{Z}^d}(n))
\label{lower2}
\endd
if and only if
\begg
\sum_{n=1}^{\infty} (a(2^n))^{Kd-d-2}<\infty.
\label{condone}
\endd
\end{thm}

The above theorems show that the behavior of distance of random walks and 
that of the distance of  Wiener  processes are very similar. This is not the 
case in two dimensions with two walkers. As it was mentioned in the 
Introduction, on $\mathbb{Z}^2$ two independent walkers will collide 
infinitely often almost surely. On the other hand, two independent standard 
Wiener processes won't collide infinitely often, their distance for $t>0$ big 
enough, will be at least $t^{-(\log t)^{\ep}}$ for any $\ep>0,$  (see 
R\'ev\'esz \cite{R13}, page 208, Remark 3). 

\smallskip
\noindent{\bf Proof of Theorem 3.1.}
It suffices to prove (\ref{iterated1}) for Wiener processes. The random walk
case (\ref{iterated2}) follows from the strong approximation in Theorem C.

Consider the case $K=2$ first. Then
$$
D_2^{\mathbb{R}^d}(t)=\|{\bf W}_1(t)-{\bf W}_2(t)\|=\sqrt{2}\|{\bf W}^*(t)\|,
$$
where ${\bf W}^*$ is a standard $d$-dimensional Wiener process. Hence by
Theorem A, (\ref{iterated1}) is true in this case.
$D_K^{\mathbb{R}^d}(t)$ is the maximum of $K\choose 2$ distances for each of
which (\ref{iterated1}) holds. This implies the upper part of the conclusion.
The lower part is immediate, namely
$$D_K^{\mathbb{R}^d}(t)\geq D_2^{\mathbb{R}^d}(t).$$
$\Box$

\noindent
{\bf Proof of Theorem 3.2.} We first  prove  the convergent part for $d=1$.
Define
\begg W_i^*(t):=\frac{\sum_{j=1}^i( W_j(t)-W_{i+1}(t))}{\sqrt{i^2+i}},
\,\,i=1,2,..., K-1. \label{insidee}
\endd
It is an easy calculation to show that $W_i^*(t) \,\, i=1,2,...,K-1 $
are independent standard Wiener processes. Hence for the $K-1$ dimensional
standard Wiener process defined by
$${\bf W}^*(t):=(W_1^*(t),W_2^*(t),...,W_{K-1}^*(t)),$$
we can  apply Theorem B with $d=K-1$ to get that, if $\{a(n),\, n=1,2,...\}$
satisfy (\ref{condone}) then for $t>t_0(\omega)$ we have
 $${\|\bf W}^*(t)\|\geq  \sqrt{t}a(t).$$
 This implies that for some $i\leq K-1$
$$\left| W_i^*(t)\right|\geq \frac{\sqrt{t}{(a(t))}}{\sqrt{K-1}},$$
which, in turn, implies that for the absolute value of one of the summands
$W_j(t)-W_{i+1}(t) , \,\,j=1,2,..., i$  of $W_i^*(t),$ we have that
$$\left|\frac{W_j(t)-W_{i+1}(t)}{\sqrt{i^2+i}}\right|\geq
\frac{\sqrt{t}{(a(t))}}{i\sqrt{K-1}},$$
implying that there exists a pair $1\leq i<j\leq K$ such that
$$ |W_i(t)-W_j(t)|\geq \frac{1}{\sqrt{K-1}} \sqrt{t}a(t).$$

\noindent
Thus, using Remark 2.1, we proved the convergent part of the theorem for $d=1$.

\smallskip
Divergent part: Suppose that
\begg
\sum_{n=1}^{\infty} (a(2^n))^{K-3}=\infty. \label{condtwo}
\endd
Then, according to Theorem B again, there is random sequence
$t_k\to \infty$ such that
$${\|\bf W}^*(t_k)\|\leq  \sqrt{t_k}a(t_k)$$
almost surely.
Thus, for all $i=1,2,...,K-1$,
$$|W_i^*(t_k)|\leq  \sqrt{t_k}a(t_k)$$
almost surely as well. Now, applying Lemma D,
we can conclude that each
 $$\{W_i(t)-W_j(t), \,i<j, \,\,i=1,2,....,K-1, \,\, j=2,3,...,K\}  \,
\label{inside2}
 $$
can be expressed as a linear combination of $W_i^*(t),\, i=1,2,... ,K-1.$ This 
in turn implies that  for all $i<j, \,\,i=1,2,....,K-1.\,\,
j=2,3,...,K,$ we have almost surely for our random sequence $t_k\to \infty$ 
that, for some appropriate constant $c_K$ we have
$$ |W_i(t_k)-W_j(t_k)|\leq c_K \sqrt{t_k}a(t_k),$$
implying that, almost surely,
$$D_K^{\mathbb{R}}(t_k)\leq c_K \sqrt{t_k}a(t_k)$$
as well, proving Theorem 3.3 in case $d=1$.

For $d>1$, we perform the previous transformation  for each coordinate  
separately. So we get $d$ times $K-1$ independent Wiener processes. For 
these $d(K-1)$  independent Wiener processes we apply Theorem B again.
Repeating the arguments in the case $d=1,$ we can finally obtain
(\ref{lower1}). Now (\ref{lower2}) follows from the strong invariance in
Theorem C.\, \,  $\Box$

\begin{remark}
We can get the limiting distribution of $D_K^{\mathbb{Z}}(n)$ via the
exact calculation  for $D_K^{\mathbb{R}}(t).$  Let
$W_i(t)/\sqrt t:=N_i$ for $ i=1,2,3,...,K$, and, conditioning on the largest
of them, say $N_1$, we  have
$$
 P\left(\frac{D_K^{\mathbb{R}}(t)}{\sqrt{t}}<z\right)
=K\int_{-\infty}^{\infty} P\left(\frac{D_K^{\mathbb{R}}(t)}{\sqrt{t}}
\leq z,\,N_i<x,\, i=2,3,...,K|N_1=x\right) P(N_1=x)\,dx$$
$$=K\int_{-\infty}^{\infty}P(x-z\leq N_i\leq x, \,i=2,3,...,K)P(N_1=x)\,dx$$
$$=K\int_{-\infty}^{\infty}\left(\Phi(x)-\Phi(x-z)\right)^{K-1} \phi(x)\,dx
=K\int_{-\infty}^{\infty}\left(\int _{x-z}^x \phi(u)\, du\right)^{K-1}
\phi(x)\,dx.
$$
where $\Phi(\cdot) $ and $\phi(\cdot)$ are the distribution and density 
functions of the standard normal random variable. Consequently, we conclude 
the following result:

\begg\lim_{n\to \infty} P\left(\frac{D_K^{\mathbb{Z}}(n)}{\sqrt n}<z\right)
=K\int_{-\infty}^{\infty}\left(\int _{x-z}^x \phi(u)\, du\right)^{K-1}
\phi(x)\,dx. \label{limitone}
\endd
It might be of interest to get the $d$-dimensional analog of this result.
\end{remark}
\section{Distance on the comb}
\renewcommand{\thesection}{\arabic{section}} \setcounter{equation}{0}
\setcounter{thm}{0} \setcounter{lemma}{0}
\bigskip

As it was mentioned in the Introduction, Krishnapur and Peres \cite{KP}
introduced a fascinating class of graphs where simple random walks continue
to be recurrent, but the respective paths of two independent random walks
meet only finitely many times with probability 1. In particular, the
2-dimensional comb lattice has this property. So, for $K$ independent walks
 $\{{\bf C}^{(i)}(n)=(C_1^{(i)}(n), C_2^{(i)}(n)) \quad i=1,2,...,K \} ,$
we want to investigate

$$D_K^{\mathbb{C}^2}(n)=\max_{i\neq j,\,\, i,j\leq K}
d({\bf C}^{(i)}(n),{\bf C}^{(j)}(n) ),$$
the maximal distance between the $K$ walkers at time $n$, where $d(x,y)$
was defined in (\ref{dist}).

The upper class result is an easy consequence of our strong approximation
in Theorem K.
\begin{thm}  For the distance of K walkers on the comb we have
\begg
\limsup_{n\to \infty} \frac{D_K^{\mathbb{C}^2}(n)}{2\sqrt{n\log\log n}}
=1 \quad {\rm a.s.}
\endd
\end{thm}
{\bf Proof.}  First we prove the theorem for two walkers. Observe that
for the second coordinates of our two walkers we have from Theorem K that
\begin{eqnarray} |C_2^{(1)}(n)-C_2^{(2)}(n)|&=&| W^*(n)-W^{**}(n)|
+O(n^{3/8+\ep})
\quad {\rm a.s.}, \nonumber \\
|C_2^{(1)}(n)+C_2^{(2)}(n)|&=&| W^*(n)+W^{**}(n)|+O(n^{3/8+\ep})
\quad {\rm a.s.},
\label{second}
\end{eqnarray}
where $W^*(n)$ and $W^{**}(n)$ are two independent standard Wiener processes.
Then both

$\displaystyle{\frac{W^*(n)-W^{**}(n)}{\sqrt{2}}}\,\,$  and
$\,\,\displaystyle{\frac{W^*(n)+W^{**}(n)}{\sqrt{2}}}\,\,$ are standard Wiener
processes again, for  which the LIL holds. Combining this with (\ref{second}),
we get that

\begg
\limsup_{n\to \infty} \frac{|C_2^{(1)}(n)-C_2^{(2)}(n)|}
{2\sqrt{n\log\log n}}=\limsup_{n\to \infty} \frac{|C_2^{(1)}(n)+C_2^{(2)}(n)|}
{2\sqrt{n\log\log n}}=1 \quad {\rm a.s.}\,.
\endd

\noindent
Applying now Corollary 2.1, (\ref{itlog1}) implies that
$|C_1^{(1)}(n)-C_1^{(2)}(n)|$
can't have a significant contribution to
$\limsup_n d({\bf C}^{(1)}(n),{\bf C}^{(2)}(n)).$
Thus the distance of the two walkers is essentially  the difference or the
sum of their  second coordinates, depending on whether they are on the same
tooth or not. Consequently, we get
\begg \limsup_{n\to \infty} \frac{D_2^{\mathbb{C}^2}(n)}
{2\sqrt{n\log\log n}}=1 \quad {\rm a.s.}\,,  \label{almost}
\endd
\noindent
so we have Theorem 4.1 for $K=2.$  From here on the proof  for $K>2$ is
exactly the same as in Theorem 3.1; by definition, $D_K^{\mathbb{C}^2}(n)$
is the maximum of $K\choose 2$ distances for each of which (\ref{almost})
holds. This implies the upper part of the theorem. The lower part is immediate,
namely $$D_K^{\mathbb{C}^2}(n)\geq D_2^{\mathbb{C}^2}(n).$$
$\Box$

 We now turn to the lower class results.
\begin{thm} For $K=2$ and any $\ep>0$
\begg
P\left(D_2^{\mathbb{C}^2}(n)\leq (1+\ep)
\frac{2^{9/4}}{3^{3/4}}n^{1/4}({\log \log n})^{3/4} \,\,{\rm i.o.}\right)=1.
\endd
\end{thm}
{\bf Proof.} The main idea of the proof is the following. Consider the second
coordinates of the two walkers. They behave like simple symmetric walks,
except that sometimes, when horizontal steps occur, they don't move. But we
know that in $n$ steps the number of vertical steps is $n(1-o(1)).$ If
they actually would move like simple symmetric walks, then, as it was mentioned
in the Introduction, they would meet infinitely often at the origin. This
means that with probability one there would be infinitely many  $n_k$ when
the second coordinates of the two walkers would be zero, and hence they both 
would be on the $x$-axis. At these occasions their distance can't be more than
what Corollary 2.1 implies, i.e., thus we arrive to our conclusion as well.

Turning to the actual proof  that the two walkers are on the $x$-axis at the 
same time infinitely often with probability  1, according to Consequence 2.3, 
for $n$ big enough, we have
\begg
P(C_2(n)=0)\geq \frac{1}{2\sqrt{n}} \label{gyokn}.
\endd

\noindent
Consider now two independent walkers
${\bf C}^{(i)}(n)=(C_1^{(i)}(n), C_2^{(i)}(n)),\,\, i=1,2$ on the comb.
Let $U$ denote the number of collisions  at zero of  their second
coordinates $C_2^{(1)}(n)$ and $C_2^{(2)}(n).$ Then by (\ref{gyokn})
$$ E(U)=E\left( \sum_{n=1}^ {\infty} I\{C_2^{(1)}(n)=C_2^{(2)}(n)=0\} \right)
\geq  \sum_{n=1}^ {\infty} \left(\frac{1}{2\sqrt{n}}\right)^2=+\infty.$$
As the number of collisions at zero of the second coordinates follows a
geometric distribution, having infinite expectation implies that there is an
infinite number of such collisions at zero. Thus, almost surely, there is a
random sequence $n_k\to \infty$ such that  $C_2^{(1)}(n_k)$ and
$C_2^{(2)}(n_k)$ are simultaneously on the backbone ($x$-axis) of the comb.
This, in turn, implies our theorem by (\ref{itlog1}) in Corollary 2.1.
$\Box$
\smallskip

As to the lower lower class (LLC)  result, first  we prove  the following 
result for $K=2.$
\begin{thm} For every $\ep>0$ for n big enough
$$ D_2^{\mathbb{C}^2}(n)> n^{1/4-\ep} \quad {\rm a.s.}$$
\end{thm}
{\bf Proof.} Define the events
\begin{eqnarray}
A_n&=& \{ D_2^{\mathbb{C}^2}(n)\leq n^{1/4-\ep}, \,\,C_1^{(1)}(n)\neq
C_1^{(2)}(n)\}, \nonumber\\
B_n&=& \{ D_2^{\mathbb{C}^2}(n)\leq n^{1/4-\ep},\,\,
C_1^{(1)}(n)= C_1^{(2)}(n) \}. \label{ab}
\end{eqnarray}
Then
$$ P(D_2^{\mathbb{C}^2}(n)\leq n^{1/4-\ep})= P(A_n)+P(B_n).$$
We show that
$$
P(A_n\, i.o.)=P(B_n\, i.o.)=0.
$$

First we give an upper bound for $P(A_n).$ To this end, we need a couple of 
lemmas.

To begin with, consider only one walk ${\bf C}(n).$  Recall the construction
of the comb walk in Section 2, where we defined $G_1,G_2,\ldots$ to be i.i.d.
geometric random variables with
$$
 P(G_1=k)=\frac1{2^{k+1}},\quad k=0,1,2,\ldots.,
$$
as the number of horizontal steps after each return to the backbone. Recall
also that $H_n$  and $V_n$ are the number of horizontal and vertical steps,
respectively, in the first $n$ steps of ${\bf C}(\cdot).$ Then it is easy to 
see that

$$H_n\leq\sum_{i=1}^{\xi_2(0,V_n)} G_i \leq \sum_{i=1}^{\xi_2(0,n)} G_i,$$
where $\xi_2(0,\cdot)$ is the local time at zero of the simple symmetric
walk $S_2(\cdot)$ of the vertical steps. Let
$$ M(C_1,n):=\max_{0\leq k\leq n} |C_1(k)|,$$
the absolute maximum of the horizontal coordinate of ${\bf C}(\cdot)$ in
$n$ steps.
\begin{lemma} For $n$ big enough
$$P(M(C_1,n) \geq n^{1/4}\log n)\leq \frac{3}{n}.
$$
\end{lemma}

\smallskip \noindent
{\bf Proof.} First we give an estimate for the upper tail of $H_n.$
Observe that, for $n$ big enough, on account of Theorem G,
we have with an appropriate constant $c>0$ and arbitrary
$\varepsilon\in (0,1/2)$ that
\begg
  P\left(\xi_2(0,n)\geq 2 \sqrt{n \log n}\right)\leq
\frac{c}{n^{2(1-\varepsilon)}}\leq \frac{1}{n}.\label{csfo}
\endd
By (\ref{csfo}) and applying Consequence 2.1,  we get that
\begin{eqnarray}
P\left(H_n\geq 3\sqrt{n \log n}\right)&\leq& P\left(\sum_{i=1}^{\xi_2(0,n)} G_i
\geq 3\sqrt{n \log n}\right) \nonumber \\ &\leq&
P\left(\sum_{i=1}^{\xi_2(0,n)} G_i\geq 3\sqrt{n \log n},\,\xi_2(0,n)
<2\sqrt{n\log n}\right)
+\frac{1}{n} \nonumber \\
&\leq& \frac{1}{n}+P\left(\sum_{i=1}^{2 \sqrt{n\log n}}G_i\geq
3\sqrt{n \log n}\right)
\leq \frac{2}{n} \label{muszaj}
\end{eqnarray}
if $n$ is big enough.

Let
$$
M(n):=\max_{0\leq k\leq n} |S(k)|
$$
\noindent
be the absolute maximum of a simple symmetric random walk $S(\cdot)$ in $n$
steps. Recall  that $C_1(n)=S_1(H_n).$  Then the well-known large deviation
result for the maximum  (see e.g.  R\'ev\'esz \cite{R13}, p. 21) and
(\ref{muszaj}) imply that
\begin{eqnarray}
P(M(C_1,n)\geq n^{1/4}\log n) &=& P(\max_{0\leq  i \leq n}|C_1(i)|
\geq n^{1/4}\log n) \nonumber\\
&\leq& \frac{2}{n} +
P(M(3\sqrt{n \log n})\geq n^{1/4}\log n)\leq \frac{3}{n}  \label{vegre}
\end{eqnarray}
if $n$ is big enough. $\Box$

\begin{lemma} Let $\mathbf{C}(n)=(C_1(n),C_2(n))$ be a random walk on
${\mathbb C}^2.$
There exists a constant $c>0$ such that we have
$$
P(C_1(n)=x,C_2(n)=y)\leq \frac{c}{n^{3/4}}\,\,\,\, {\it for\,\, all}\,\,
(x,y) \,\,{\it with}\,\, |y|\leq n^{1/2-\ep}.
$$
\end{lemma}

\smallskip \noindent
{\bf Proof.} In what follows, unimportant constants will be denoted by $c,$
whose value might change from line to line. For simplicity, we work with
even coordinates, and that, as we remarked earlier, does not restrict
generality. Recall that $H_n$ is the number of horizontal steps in the first
$n$ steps of the walk.
We have

\begin{eqnarray} P(C_1(2n)=2r, C_2(2n)=2j)
&=&\sum_k P(C_1(2n)=2r, C_2(2n)=2j|H_{2n}=2k)P(H_{2n}=2k)\nonumber \\
&=&\sum_k P(S_1(2k)=2r)P(C_2(2n)=2j|H_{2n}=2k)P(H_{2n}=2k)\nonumber\\
&\leq& \sum_k P(S_1(2k)=0)P(C_2(2n)=2j|H_{2n}=2k)P(H_{2n}=2k)\nonumber\\
&=&P(C_1(2n)=0, C_2(2n)=2j).
\end{eqnarray}

The second equality above follows from the fact that when the number of
horizontal steps are fixed, then $C_1(\cdot)$ is a simple symmetric walk,
denoted by $S_1(\cdot)$, which is independent of the second coordinate. The
above inequality, on the other hand, is true, as
$$\max_{-k\leq r\leq k} P(S_1(2k)=2r)=P(S_1(2k)=0).$$
To finish the proof, observe that by Lemma H, for $j\neq 0,$
$$P(C_1(2n)=0, C_2(2n)=2j)=p^{{\mathbb C}^2}((0,0),(0,2j),2n)
=\frac{1}{2}p^{{\mathbb C}^2}((0,2j),(0,0),2n).$$
Now our lemma follows from Consequence 2.2.
$\Box$

Returning now to the proof of Theorem 4.3, we can give the following upper
bound for $P(A_n).$ Define the event
$$U_n:=\{M(C^{(i)}_1,n)\leq n^{1/4}\log n,\,\,i=1,2\}.$$

\noindent
Then, by Lemma 4.1, we have for  $U_n^c,$ the complement of $U_n,$ that

$$P(U_n^c)\leq \frac{6}{n}.$$
Define now the set of pairs of points on the comb
$$Q(n,\ep)=\{(x_1,y_1),(x_2,y_2): |x_i|\leq n^{1/4}\log n,\,
i=1,2,\,\, |x_1-x_2|\leq n^{1/4-\ep},\,
$$
$$
|y_1|\leq n^{1/4-\ep},\,
|y_2|\leq n^{1/4-\ep}\}.$$
By Lemma 4.2
\begin{eqnarray}
 P(A_n)&\leq& P(|C_2^{(1)}(n)|\leq n^{1/4-\ep},\,\, |C_2^{(2)}(n)|
\leq n^{1/4-\ep},\,\, |C_1^{(1)}(n)-C_1^{(2)}(n)|\leq n^{1/4-\ep}) \nonumber \\
 &\leq& P(|C_2^{(1)}(n)|\leq n^{1/4-\ep},\,\, |C_2^{(2)}(n)|
\leq n^{1/4-\ep},\,\, |C_1^{(1)}(n)-C_1^{(2)}(n)|\leq n^{1/4-\ep} ,
U_n)+\frac{6}{n}\nonumber \\
 &\leq& \sum_{Q(n,\ep)} P(\mathbf{C}^{(1)}(n)=(x_1,y_1),
\mathbf{C}^{(2)}(n)=(x_2,y_2))+\frac{6}{n}
 \nonumber \\
 &\leq&
 16(n^{1/4}\log n) (n^{1/4-\ep})^3\left(\frac{c}{n^{3/4}}\right)^2+
\frac{6}{n}\leq \frac{c\log n}{n^{1/2+3\ep}}, \label{ime}
 \end{eqnarray}
\noindent
implying that for the subsequence $n_k=k^{\alpha}$,
with any $\alpha>2$,
$$\sum_{k=1}^{\infty} P(A_{n_k})< \infty .$$

\noindent
This implies that, almost surely for $k\geq k_0(\omega),$
$A_{n_k}$ does not occur, which in turn means that if the two walkers are on
different teeth, then either
$$A_1(k):=\{|C_2^{(1)}(n_k)|\geq n_k^{1/4-\ep}\},\,\, \,\,{\rm or}\,\,
A_2(k):=\{|C_2^{(2)}(n_k)|\geq n_k^{1/4-\ep}\},\,\,\,\,
$$
$$
{\rm or}\,\,A_3(k):=\{|C_1^1(n_k)-C_2^1(n_k)|\geq n_k^{1/4-\ep}\}$$
will occur. We want to show that we can select $\alpha>2$ such that for any
$n,$ with $n_k\leq n\leq n_{k+1},$  if one of the events
$\{A_i(k)\,\, i=1,2,3\}$ occurs, then
$$D_2^{\mathbb{C}^2}(n)\geq n^{1/4-\ep}$$
will occur as well, as long as two walkers are on different teeth. Since 
$n_k=k^{\alpha}$, we have $n_{k+1}-n_k\sim \alpha k^{\alpha-1}.$ So we have 
to show that in $\alpha k^{\alpha-1}$ steps the increments of the three 
processes in the events $\{A_i(k)\,\, i=1,2,3\}$ are less than $n_k^{1/4-\ep}.$
The first two of these three processes are simple symmetric walks, while the
third one is  the difference of two simple symmetric walks, but with a much
smaller number of steps (as there are  possible vertical excursions when the
horizontal move pauses). By Theorem E the increment of these walks in
$\alpha k^{\alpha-1}$ steps is almost surely less than

$$k^{\frac{\alpha-1}{2}}\log k,$$
 while
 $$n_k^{1/4-\ep}=k^{(1/4-\ep)\alpha}.$$
So we need to have
$$\frac{\alpha-1}{2}< \alpha\left(\frac{1}{4}-\ep\right),$$
which is equivalent to $\alpha<\frac{2}{1+4\ep}.$
On the other hand, for the convergence of $\sum_kP(A_{n_k})$ we need that
$\alpha(1/2+3\ep)\geq 1$ should hold, which is equivalent to
$\alpha>\frac{2}{1+6\ep}.$
So, for any $\ep>0$, we can find an appropriate
$$\frac{2}{1+6\ep}<\alpha <\frac{2}{1+4\ep} ,$$
 and  conclude by the Borel-Cantelli Lemma that
 \begg
 P(A_n \,i.o.)=0. \label{finite}
\endd

To show that $P(B_n\, i.o.)=0$,  recall the definition of $B_n$ in
(\ref{ab}). Equivalently,
$$B_n=\{C_1^{(1)}(n)=C_1^{(2)}(n),\, |C_2^{(1)}(n)-C_2^{(2)}(n)|\leq
n^{1/4-\varepsilon}\}=B_n^{1}\cup B_n^{2},
$$
where
$$
B_n^{1}=\{C_1^{(1)}(n)=C_1^{(2)}(n),\, |C_2^{(1)}(n)-C_2^{(2)}(n)|\leq
n^{1/4-\varepsilon}, \max(|C_2^{(1)}(n)|, |C_2^{(2)}(n)|)\leq n^{1/4}\},
$$
$$
B_n^{2}=\{C_1^{(1)}(n)=C_1^{(2)}(n),\, |C_2^{(1)}(n)-C_2^{(2)}(n)|\leq
n^{1/4-\varepsilon}, \max(|C_2^{(1)}(n)|, |C_2^{(2)}(n)|)> n^{1/4}\}.
$$

\noindent
We first show  that $P(B_n^{1}\, i.o.)=0.$
$P(B_n^{1})$ can be estimated similarly to $P(A_n)$ in (\ref{ime}).
We obtain
$$
P(B_n^{1})\leq \frac{c\log n}{n^{3/4}}.
$$
Choosing $n_k=k^\alpha$ with
$$\frac{4}{3}<\alpha<\frac{2}{1+4\varepsilon},
$$ we have $P(B_{n_k}^{1}\, i.o.)=0,$
i.e., there exists a $k_0$ such that $B_{n_k}^{1}$ does not occur almost surely
if $k\geq k_0$. As we already proved that  $P(A_{n}\, i.o.)=0,$ and we
will prove that  $P(B_{n}^{2}\, i.o.)=0,$ we may assume that for
$k\geq k_0$ neither $A_{n_k}$ nor $B_{n_k}^2$ occur.
Consequently, it suffices to consider the case when
$$
|C_2^{(1)}(n_k)-C_2^{(2)}(n_k)|>n_k^{1/4-\varepsilon},
$$
since otherwise, either
$$
C_1^{(1)}(n_k)\neq C_1^{(2)}(n_k),$$ in which case
$A_{n_k}$ does not occur, or
$$
\max(|C_2^{(1)}(n_k)|,|C_2^{(2)}(n_k)|)>n_k^{1/4},
\qquad C_1^{(1)}(n_k)=C_1^{(2)}(n_k),
$$
in which case $B_{n_k}^2$ does not occur.
Now let $n_k\leq n<n_{k+1}$. We have to show that if
$C_1^{(1)}(n)=C_1^{(2)}(n)$ and
$\max(|C_2^{(1)}(n)|, |C_2^{(2)}(n)|)\leq n^{1/4}$, then
\begg
|C_2^{(1)}(n)-C_2^{(2)}(n)|>n^{1/4-\varepsilon},
\label{big}
\endd
i.e., $B_n^1$ does not occur with probability 1 for large $n$.
The increments of $C_2(\cdot)$ in $(n_k,n_{k+1})$ are almost surely less than
$$
k^{\frac{\alpha-1}{2}}\log k<n_k^{1/4-\varepsilon},
$$
so it can be seen that (\ref{big}) holds, i.e., $B_n^1$ does not occur,
so $P(B_n^1\, i.o.)=0$.

To prove $P(B_n^2\, \,i.o.)=0,$ we need the following Lemma.
\begin{lemma}
Let $E_n$ and $B_n$ be two sequences of events on the same probability space.
Introduce the notations
\begg
\quad E_n^*:=\bigcup_{j=n}^\infty E_j. \quad {\rm and}
\quad B_{n,m}^*=B_n\cap B_{n-1}^c\cap\ldots\cap B_m^c,\quad m<n,\quad
{\rm with}\quad B_{n,n}^*=B_n,\label{cond2}
\endd
where $B^c$ denotes the complement of $B$.
Assume that
\begg
P(E_n\, i.o.)=0,  \label{cond1}
\endd
and
\begg
P(E_n^*|B_{n,m}^*)\geq C>0\label{cond3}
\endd
for large enough $m\leq n$ with some constant $C$. Then we also have
$$
P(B_n\, i.o.)=0.
$$
\end{lemma}
{\bf Proof.} It is known that $P(E_n\, i.o.)=0$ is equivalent to
$\lim_{n\to\infty}P(E_n^*)=0$.
 $B_{n,m}^*$ means that $n$ is the
first index, when $B_i$ occurs with $i\geq m$. Then by (\ref{cond3}) we have 
that
$$
P(E_n^*\cap B_{n,m}^*)=P(E_n^*|B_{n,m}^*)P(B_{n,m}^*)\geq CP(B_{n,m}^*).
$$
Since $B_{k,m}^*$ are disjoint for different $k$, and $B_k\supset B_{k,m}^*$,
we have
\begin{eqnarray}
P(\bigcup_{k=m}^{\infty} E_k^*\cap B_k)&\geq&
P(\bigcup_{k=m}^\infty E_k^*\cap B_{k,m}^*)
=\sum_{k=m}^\infty P(E_k^*\cap B_{k,m}^*)\nonumber  \\
&\geq& C\sum_{k=m}^\infty P(B_{k,m}^*)=CP(\bigcup_{k=m}^\infty B_{k,m}^*)
=CP(\bigcup_{k=m}^\infty B_k). \label{harom}
\end{eqnarray}
By  (\ref{cond1}), $P(E_m^*\cap B_m\, i.o.)=0$ as well, or equivalently,
$$
\lim_{m\to\infty}P(\bigcup_{k=m}^\infty E_k^*\cap B_k)=0.
$$
Consequently, by (\ref{harom}),
$$
\lim_{m\to\infty}P\left(\bigcup_{k=m}^\infty B_k\right)=0,
$$
hence $P(B_n\, i.o.)=0.$  $\Box$

To complete the proof of Theorem 4.3, we have to prove
$P(B_n^2\,\, i.o.)=0$. Recall the result of Krishnapur and Peres \cite{KP}
that $P({\bf C}^{(1)}(n)={\bf C}^{(2)}(n)\, \,i.o.)=0$. Similarly, it can be
shown that for
\begg
E_n=\{C_1^{(1)}(n)=C_1^{(2)}(n),\, |C_2^{(1)}(n)-C_2^{(2)}(n)|\leq 1\},
\label{newe}
\endd
we have also $P(E_n\, \,i.o.)=0$. To apply Lemma 4.3 with $E_n$ as in
(\ref{newe}) and $B_n$ replaced by $B_n^{2}$, we have to prove that
\begg
P(E_n^*|B_{n,m}^{2*})\geq c
\label{pos}
\endd
with some constant $c>0$, by showing that if ${\bf C}^{(1)}(n)$ and
${\bf C}^{(2)}(n)$ are as in $B_n^{2}$, then before returning to the
backbone, with positive probability they either meet at some point, or are
at distance 1. Now define
\begg
\tau_1=\min\{k\geq 0:\, |C_2^{(1)}(n+k)-C_2^{(2)}(n+k)|\leq 1\},
\label{tau1}
\endd
\begg
\tau_2^{(j)}=\min\{k\geq 0:\, C_2^{(j)}(n+k)=0\},\quad j=1,2,\quad
\tau_2=\min(\tau_2^{(1)},\, \tau_2^{(2)}).
\label{tau2}
\endd
Since, under the condition $B_{n,m}^{2*},$ both $C_2^{(1)}(n)$ and
$C_2^{(2)}(n)$ are either positive or negative, we have
$$
P(E_n^*|B_{n,m}^{2*}) \geq P(\tau_1<\tau_2|B_{n,m}^{2*}).
$$
What we have to show is that this  last probability can be bounded from below
by a positive constant. This will be achieved by estimating the distributions
of $\tau_1$ and $\tau_2$, and applying Theorem L in Section 2, since these
distributions are equivalent in terms of $\beta(r)$ as in Theorem L.
\begin{lemma}
$$
\lim_{m\to\infty}P(\tau_1<n^{1/2-\varepsilon}|B_{n,m}^{2*})=1,
\qquad \lim_{m\to\infty}P(\tau_2>n^{1/2-\varepsilon/2}|B_{n,m}^{2*})=1.
$$
\end{lemma}
{\bf Proof.}
\begin{eqnarray*}
&P&(\tau_1<n^{1/2-\varepsilon}|B_{n,m}^{2*}) \nonumber  \\
&=&\sum_{z_1,z_2}P(\tau_1<n^{1/2-\varepsilon}|B_{n,m}^{2*},\, C_2^{(1)}(n)=z_1,
\, C_2^{(2)}(n)=z_2)P(C_2^{(1)}(n)=z_1,\, C_2^{(2)}(n)=z_2|B_{n,m}^{2*}),
\end{eqnarray*}
and
\begin{eqnarray*}
&P&(\tau_2>n^{1/2-\varepsilon/2}|B_{n,m}^{2*})=\nonumber  \\
&=&\sum_{z_1,z_2}P(\tau_2>n^{1/2-\varepsilon/2}\,|B_{n,m}^{2*},\,
C_2^{(1)}(n)=z_1,\,
C_2^{(2)}(n)=z_2)P(C_2^{(1)}(n)=z_1,\, C_2^{(2)}(n)=z_2|B_{n,m}^{2*}),
\nonumber  \\
\end{eqnarray*}
where the summation $\sum_{z_1,z_2}$ stands for all permissible values of
$z_1,z_2$, under the condition $B_{n,m}^{2*}$. In fact, under the condition
$\{B_{n,m}^{2*},\, C_2^{(1)}(n)=z_1,\, C_2^{(2)}(n)=z_2\},$
$$C_2^{(1)}(n+k), \, C_2^{(2)}(n+k),\quad k=0,1,\ldots,\tau_2,$$
are two independent simple random walks, starting at $z_1$ and $z_2$,
respectively, and avoiding 0 before $\tau_2$. Moreover, since under the
above condition, both $C_2^{(1)}(n)$ and $C_2^{(2)}(n)$ are either positive
or negative,
$$
|C_2^{(1)}(n+k)-C_2^{(2)}(n+k)|
$$
for $0\leq k\leq \tau_1$ behaves also as a simple random walk with even number
of steps, starting at $|z_1-z_2|$, $\tau_1$ being the first hitting time of
zero or one, depending on the parity of $|z_1-z_2|$. The distribution of
$\tau_1$ is equivalent to that of the first hitting time of $|z_1-z_2|$ or
$|z_1-z_2|-1$ of a simple random walk, starting from $0$ and considering
even number of steps. Denoting by $S(\cdot)$ a simple random walk on the
line, it can be seen that
$$
P(\tau_1<n^{1/2-\varepsilon}|B_{n,m}^{2*})\geq
P(\tau_1<n^{1/2-\varepsilon}|C_1^{(1)}(n)=C_1^{(2)}(n),\,
|C_2^{(1)}(n)-C_2^{(2)}(n)|=2[n^{1/4-\varepsilon}])
$$
$$
=P(\min\{k\geq 0: S(2k)=0\}<n^{1/2-\varepsilon}|S(0)=2[n^{1/4-\varepsilon}])=
P(\beta(2[n^{1/4-\varepsilon}])<2n^{1/2-\varepsilon}),
$$
since $2\min\{k\geq 0: S(2k)=0\}$ under the condition
$S(0)=2[n^{1/4-\varepsilon}]$ has the same distribution as
$\beta(2[n^{1/4-\varepsilon}])$ in Theorem L.

Concerning $\tau_2$, suppose that  $|C_2^{(i)}(n)|<
|C_2^{(j)}(n)|. $  Then it suffices to consider the time when the random walk
$C_2^{(i)}(n+k),\, k\geq 0,$ reaches $0$ , otherwise the two random walks 
will meet before $\tau_2$, i.e.,  $\tau_1<\tau_2$. Under the condition 
$B_{n,m}^{2*}$, we have
$$
P(\tau_2>n^{1/2-\varepsilon/2}|B_{n,m}^{2*})\geq P(\min\{k\geq 0:\,
S(k)=0\}>n^{1/2-\varepsilon/2}|S(0)=[n^{1/4}])
$$
$$
=P(\beta([n^{1/4}])>n^{1/2-\varepsilon/2}),
$$
where $\beta(\cdot)$ is defined in (\ref{beta}) and $[\cdot]$ denotes
integral part. Now Lemma 4.4 follows from the limiting distributions of
hitting times in Theorem L of Section 2. $\ \ \ \ \ \Box$

To complete the proof of Theorem 4.3, it follows from Lemma 4.4 that
$$
P(\tau_1<\tau_2|B_{n,m}^{2*})\geq P(\tau_1<n^{1/2-\varepsilon},\,
\tau_2>n^{1/2-\varepsilon/2}|B_{n,m}^{2*})
$$
$$
\geq P(\tau_1<n^{1/2-\varepsilon}|B_{n,m}^{2*})
-P(\tau_2<n^{1/2-\varepsilon/2}|B_{n,m}^{2*})\geq c>0,
$$
where the first term tends to 1, the second term tends to zero, as
$n\to\infty$, so the difference is greater than a positive constant $c$.
This means that the two random walks will meet after time
$n$ with positive probability, i.e., (\ref{pos}) holds. Hence, using
Lemma 4.3, we have $P(B_n^2\, i.o.)=0.$ This completes the proof of
Theorem 4.3. $\Box$

Concerning the lower classes for more than $2$ walkers, we have the following
result.
\begin{thm}
Let $a(n)$ be a nonincreasing nonnegative function. Then, for $K\geq 3$
\begg
\sqrt{n} a(n) \in {\rm LLC}(D_K^{\mathbb{C}^2}(n))
\endd
if and only if
\begg
\sum_{n=1}^{\infty} (a(2^n))^{K-2}<\infty. \label{condone2}
\endd
\end{thm}
{\bf Proof.} First assume that
$$
\sum_{n=1}^{\infty} (a(2^n))^{K-2}=\infty.
$$
Then, by the Dvoretzky-Erd\H os Theorem (Theorem B in Section 2) and the strong
approximation in Theorem K, we have
$$
\max_{1\leq j\leq K}|C_2^{(j)}(n)|\leq
\sqrt{\sum_{j=1}^K (C_2^{(j)}(n))^2}\leq n^{1/2}a(n)/3
$$
infinitely often with probability 1, since we have also
$$
\sum_{n=1}^{\infty} (a(2^n)/3)^{K-2}=\infty.
$$
By (\ref{itlog1}) of Corollary 2.1 for the horizontal distance, we have
$$
|C_1^{(i)}(n)-C_1^{(j)}(n)|\leq n^{1/2}a(n)/3
$$
for all $1\leq i,j\leq K$ and all large enough $n,$ with probability 1.
Then
$$
d({\bf C}^{(i)}(n),\, {\bf C}^{(j)}(n))\leq
|C_2^{(i)}(n)|+|C_2^{(j)}(n)|+|C_1^{(i)}(n)-C_1^{(j)}(n)|,
$$
and, consequently, we have
$$
D_K^{\mathbb{C}^2}(n)\leq 2\max_{1\leq j\leq K}|C_2^{(j)}(n)|
+\max_{1\leq i,j\leq K}|C_1^{(i)}(n)-C_1^{(j)}(n)|\leq n^{1/2}a(n)
$$
infinitely often with probability 1. This verifies the first part of
Theorem 4.4.

To show the other part, i.e., assuming that
$$
\sum_{n=1}^{\infty} (a(2^n))^{K-2}<\infty,
$$
we have to prove
\begg
D_K^{\mathbb{C}^2}(n)> n^{1/2}a(n)
\label{d3}
\endd
for all large enough $n$ with probability 1. The idea is similar to the
proof of Theorem 4.3 concerning the event $B_n^{2}$ in the case when the 2
random walks are on the same tooth at time $n$. In fact, we show that one
of the random walks has to be high on some tooth by the Dvoretzky-Erd\H os
Theorem B, and no other random walks can be close to this one on the same
tooth. We note that the constants in the following proof are not too important,
one could also choose different suitable  constants.

Assume that we have $K$ independent random walks
${\bf C}^{(i)}(\cdot),\, i=1,\ldots,K$ on the comb. By Theorem B and
Theorem K, we have
$$
\max_{1\leq j\leq K}(|C_2^{(j)}(n)|)>5n^{1/2}a(n).
$$
If there is no other random walk on the same tooth than the one taking the
above maximum of $|C_2^{(j)}(n)|$, then obviously (\ref{d3}) holds. So we can
consider the case when two random walks are on the same tooth at time $n,$
and one of them is higher than $5n^{1/2}a(n)$. For fixed $1\leq i,j\leq K,$
define the event
$$
B_n^{i,j}=\{C_1^{(i)}(n)=C_1^{(j)}(n)\, , \, |C_2^{(i)}(n)-C_2^{(j)}(n)|\leq
n^{1/2}a(n),\, \max(|C_2^{(i)}(n)|,|C_2^{(j)}(n)|)>5n^{1/2}a(n)\}.
$$

We show that for fixed $i,j$, $P(B_n^{i,j}\, i.o.)=0$, by applying
Lemma 4.3, with $E_n$ as in (\ref{newe}). Note that $B_n^{i,j}$ implies that
$\min(|C_2^{(i)}(n)|\, ,\, |C_2^{(j)}(n)|)>3n^{1/2}a(n)$.

Define $\tau_1$ and $\tau_2$ as in (\ref{tau1}) and (\ref{tau2}), with the
obvious modification that $(1),(2)$ should be $(i),(j)$. Here again, we can
consider $\tau_2$ to be the time when the lower one  of $|C_2^{(i)}(n)|$ and
$|C_2^{(j)}(n)|$ reaches zero. Then, similarly to the proof of Theorem 4.3,
using Theorem L,
$$
P(\tau_1<\tau_2|B_{n,m}^{i,j*})
\geq P(\tau_1<2na^2(n)|B_{n,m}^{i,j*})-P(\tau_2<9na^2(n)/2|B_{n,m}^{i,j*})
$$
$$
\geq P(\beta([2n^{1/2}a(n)])<4na^2(n))-P(\beta([3n^{1/2}a(n)])<9na^2(n)/2)>c
$$
with some positive constant $c$. Hence, by Lemma 4.3,
$P(B_n^{i,j}\, i.o.)=0$ for all $1\leq i,j\leq K$. This means that there is
at least one distance $|C_2^{(i)}(n)-C_2^{(j)}(n)|$ larger than
$n^{1/2}a(n)$, for all large $n$ with probability 1, so (\ref{d3}) follows.
This completes the proof of Theorem 4.4.  $\Box$

\bigskip
\noindent
{\bf Acknowledgements.} The authors are
grateful to Mikl\'os Cs\"org\H o for useful remarks.


\begin{thebibliography}{}

\bibitem{BPS}
Barlow, M. T., Peres, Y. and Sousi, P.: Collision of random
walks. \textit{Annales de l'Institut Henri Poincar\'e- Probabilit\'es
et Statisiques} \textbf{48} (2012), 922-946.

\bibitem{BE}
Bertacchi, D.: Asymptotic behaviour of the simple random
walk on the 2-dimensional comb. \textit{Electron. J. Probab.}
\textbf{11} (2006), 1184--1203.

\bibitem{BZ0}
Bertacchi, D. and Zucca, F.: Equidistribution of random
walks on spheres. \textit{J. Stat. Phys.} \textbf{94} (1999), 91--111.

\bibitem{BZ}
Bertacchi, D. and Zucca, F.: Uniform asymptotic
estimates of transition probabilities on combs. \textit{J. Aust. Math.
Soc.} \textbf{75} (2003), 325--353.

\bibitem{CCFR08}
Cs\'aki, E., Cs\"org\H o, M., F\"oldes, A. and R\'ev\'esz, P.:
Strong limit theorems for a simple random walk on the
2-dimensional comb. \textit{Electron. J. Probab.} \textbf{14}  (2009),
2371--2390.

\bibitem{CsF83}
Cs\'aki, E. and F\"oldes, A.: How big are the increments of the local time
of a recurrent random walk? \textit{Z. Wahrsch. Verw. Gebiete} \textbf{65}
(1983), 307--322.


\bibitem{CR81}
Cs\"org\H o, M. and R\'ev\'esz, P.: \textit{Strong Approximations in
Probability and Statistics.} Akad\'emiai Kiad\'o, Budapest and
Academic Press, New York, 1981.

\bibitem{DE}
Dvoretzky, A. and Erd\H{o}s, P.: Some problems on random walk in
space. In: \textit{Proc. of the Second Berkeley Sympos. Math. Statist. and
Probability}, Berkeley, Calif. 1950, pp. 353--367.

\bibitem{FE}
Feller, W.: \textit{An Introduction to Probability Theory and its
Applications, Vol. 1}, 3rd edn. John Wiley and Sons, NewYork-London-
Sidney, 1968.

\bibitem{KP}
Krishnapur, M. and Peres, Y.: Recurrent graphs where two
independent random walks collide finitely often. \textit{Electron. Comm.
Probab.} \textbf{9} (2004), 72--81.

\bibitem{P}
P\'olya, G.: \"Uber eine Aufgabe der Wahrsheinlichkeitsrechnung
betreffend die Irrfahrt in Strassennetz. \textit{Math. Ann.} \textbf{84}
(1921), 149-160.

\bibitem{R13}
R\'ev\'esz, P.: \textit{Random Walk in Random and Non-Random Environments},
3rd edn. World Scientific, Singapore, 2013.

\bibitem{TO}
T\'oth, B.: No more than three favorite sites for
simple random walk. \textit{Ann. Probab.} \textbf{29} (2001), 484--503.

\bibitem{WH}
Weiss, G. H. and Havlin, S.: Some properties of a random
walk on a comb structure. \textit{Physica A} \textbf{134}  (1986), 474--482.
\end{thebibliography}
\end{document}